\documentclass[12pt]{article}
\usepackage{amssymb}
\usepackage{amsthm}
\usepackage{amsmath}
\usepackage{graphics}
\usepackage{epsfig}
\usepackage{float}

\def\Que{{\mathbb Q}}
\def\En{{\mathbb N}}
\def\Zee{{\mathbb Z}}
\def\mod#1 #2{#1\ ({\rm mod}\ #2)}

\newtheorem{theorem}{Theorem}

\newtheorem{lemma}[theorem]{Lemma}

\newtheorem{cor}[theorem]{Corollary}
\newtheorem{corollary}[theorem]{Corollary}

\newcommand{\set}[1]{\left\{#1\right\}}

\title{Van der Waerden's Theorem and 
Avoidability in Words}

\author{Yu-Hin Au\\
{\small Department of Combinatorics \& Optimization,
University of Waterloo,}\\[-2pt]{\small
Waterloo, Ontario  N2L 3G1,
Canada}\\
{\tt \small yau@uwaterloo.ca} \\
\ \\
Aaron Robertson\\
{\small Department of Mathematics,
Colgate University,
Hamilton, NY  13346, USA} \\
{\tt \small arobertson@colgate.edu} \\
\ \\
Jeffrey Shallit\\
{\small School of Computer Science,
University of Waterloo,
Waterloo, Ontario  N2L 3G1,
Canada}\\
{\tt \small shallit@cs.uwaterloo.ca}
}

\begin{document}

\maketitle

\begin{abstract}
Pirillo and Varricchio, and independently, Halbeisen and Hungerb\"uhler
considered the following problem, open since 1994:
Does there exist an infinite word ${\bf w}$
over a finite subset of $\Zee$ such that $\bf w$ contains no two consecutive
blocks of the same length and sum?  We consider some variations
on this problem in the light of
van der Waerden's theorem on arithmetic progressions.
\end{abstract}

\section{Introduction}
\label{sec1}

Avoidability problems play a large role in combinatorics on words
(see, e.g., \cite{Lothaire:1983}).
By a {\it square} we mean a nonempty word of the form
$xx$, where $x$ is a word; an example in English is {\tt murmur}.
A classical avoidability problem is the following:
Does there exist an infinite word
over a finite alphabet that contains no squares?
It is easy to see that no such word
exists if the alphabet size is $2$ or less, but if the alphabet size
is $3$, then such a word exists, as proven by Thue \cite{Thue:1906,Thue:1912}
more than a century ago.  

     An {\it abelian square} is a nonempty word of the form $xx'$
where $|x| = |x'|$ and $x'$ is a permutation of $x$.  An example
in English is {\tt reappear}.
In 1961, Erd\H{o}s \cite{Erdos:1961}
asked: Does there exist an infinite word over a finite alphabet
containing no 
abelian squares?
Again, it is not hard to see that this
is impossible over an alphabet of size less than $4$.  Evdokimov
\cite{Evdokimov:1968} and Pleasants \cite{Pleasants:1970} gave
solutions for alphabet size $25$ and $5$, respectively, but it was not until
1992 that Ker\"anen \cite{Keranen:1992} proved that an infinite word
avoiding abelian squares does indeed
exist over a $4$-letter alphabet.

     Pirillo and Varricchio \cite{Pirillo&Varricchio:1994}, and
independently, Halbeisen and Hungerb\"uhler \cite{Halbeisen&Hungerbuhler:2000},
suggested yet another variation.  Let a {\it sum-square} be a
factor of the form $x x'$ with $|x| = |x'|$ and $\sum x = \sum x'$, where
by $\sum x$ we mean the sum of the entries of $x$.  
Is it possible to construct an infinite
word over a finite subset of $\Zee$ that contains no sum-squares?
This very interesting
question has been open for 15 years. Halbeisen and Hungerb\"uhler
observed that the answer is ``no'' if we omit the condition 
$|x| = |x'|$.  Their tool was a famous one from combinatorics:  namely,
van der Waerden's theorem on arithmetic progressions
\cite{van.der.Waerden:1927}.

\begin{theorem} \rm(van der Waerden)\it \,
      Suppose $\En$ is colored using a finite number of colors.
Then there exist arbitrarily long monochromatic
arithmetic progressions.
\end{theorem}

In this note we consider several variations on the problem of
Pirillo-Varricchio and Halbeisen-Hungerb\"uhler (the {\it PVHH
problem}, for short). In Section~\ref{sec2}, we show there is
no infinite abelian squarefree word in which the difference between the
frequencies of any two letters is bounded above by a constant.
Section~\ref{sec3} deals with the problem of avoiding sum-squares,
modulo $k$. While it is known there is no infinite word with this
property (for any $k$), we show that there is an infinite word over
$\{-1,0,1\}$ that is squarefree and avoids all sum-squares in which the
sum of the entries is non-zero.

In Section~\ref{sec4}, we provide upper and lower bounds on the length of any
word over $\mathbb{Z}$ that avoids sum-squares (and 
higher-power-equivalents) modulo $k$. We conclude with some
computational results in Section~\ref{sec5}.

\section{First Variation}
\label{sec2}

We start with an infinite word $\bf w$
already known to avoid abelian squares
(such as Ker\"anen's, or other words found by Evdokimov \cite{Evdokimov:1968}
or Pleasants \cite{Pleasants:1970}) over some finite  alphabet
$\Sigma_k = \lbrace 0, 1, \ldots, k-1 \rbrace$.
We then choose an integer base $b \geq 2$ and
replace each occurrence of $i$ in $\bf w$ with $b^i$, obtaining a new
word $\bf x$.  If there were
no ``carries'' from one power of $b$ to another, then $\bf x$ would
avoid sum-squares.  We can avoid problematic
``carries'' if and only if, whenever
$x x'$ is a factor with $|x| = |x'|$, then the number of occurrences of
each letter in $x$ and $x'$ differs by less than $b$.  In other words,
we could solve the PVHH problem
if we could find an abelian squarefree word such that
the difference in the number of occurrences between the 
most-frequently-occurring and least-frequently-occurring letters in any
prefix is bounded.  As we will see, though, this is impossible.

More generally, we consider the frequencies of letters in
abelian power-free words.
By an {\it abelian $r$-power} we mean 
a factor of the form $x_1 x_2 \cdots x_r$, where
$|x_1| = |x_2| = \cdots = |x_r|$ and
each $x_i$ is a permutation of $x_1$.  For example, the English word
{\tt deeded} is an abelian cube.

We introduce some notation.
For a finite
word $w$, we let $|w|$ be the length of $w$ and let $|w|_a$ be the number
of occurrences of the letter $a$ in $w$.  Let $\Sigma = \lbrace a_1, a_2,
\ldots, a_k \rbrace$ be a finite ordered alphabet.  Then
for $w \in \Sigma^*$, we let $\psi(w)$ denote the vector
$(|w|_{a_1}, |w|_{a_2}, \ldots, |w|_{a_k})$.  The map $\psi$ is sometimes
called the {\it Parikh map}.  For example, if $\Sigma = 
\lbrace {\tt v, l, s, e} \rbrace$, then
$\psi({\tt sleeveless}) = (1,2,3,4)$.

For a vector $u$, we let $u_i$ denote the $(i+1)^{\mathrm{st}}$ entry, so that
$u= (u_0, u_1, \ldots, u_{k-1})$.  If $u$ and $v$ are two vectors 
with real entries, we define their $L^\infty$ distance $\mu(u,v)$ to be 
$$ \max_{0 \leq i < k} |u_i - v_i|.$$

If ${\bf w} = b_1 b_2 \cdots$ is an infinite word, with each $b_i \in \Sigma$,
then by
${\bf w}[i]$ we mean the symbol $b_i$ and by ${\bf w}[i..j]$
we mean the word $b_i b_{i+1} \cdots b_j$.  Note that if $i = j+1$,
then ${\bf w}[i..j] = \epsilon$, the empty word.

\begin{theorem}\label{p1}
Let $\bf w$ be an infinite word over the finite alphabet
$\set{0, 1, \ldots, k-1}$ for some $k \geq 1$.
If there exist a vector 
$v \in \Que^k$ and a positive integer $M$ such that
\begin{equation}\label{e1}
\mu( \psi( {\bf w}[1..i]), iv) \leq M
\end{equation}
for all $i \geq 0$, 
then $\bf w$ contains
an abelian $\alpha$-power for every integer $\alpha \geq 2$.
\end{theorem}

\begin{proof}
First, note that 
\begin{equation}
\sum_{0 \leq i < k} v_i = 1.
\label{sumv}
\end{equation}
For otherwise we have $\sum_{0 \leq i < k} v_i = c \not = 1$, and then
$\mu(\psi({\bf w}[1..i]), iv)$ is at least $|c-1| {i \over k}$, and hence unbounded
as $i \rightarrow \infty$.  

For $i \geq 0$, define
$X^{(i)} = \psi( {\bf w}[1..i] ) - iv$.  Then
\begin{eqnarray}
X^{(i+j)} - X^{(i)} &=& ( \psi( {\bf w}[1..i+j] ) - (i+j) v)
	- (\psi( {\bf w}[1..i] ) - iv)  \nonumber \\
&=& \psi( {\bf w}[i+1..i+j] ) - jv \label{jv}
\end{eqnarray}
for integers $i, j \geq 0$.
For $i \geq 0$,
define $\Gamma(i)$ to be the vector with $\binom{k}{2}$ entries
given by $X^{(i)}_l - X^{(i)}_m$ for $0 \leq l < m < k$. 

From~\eqref{e1}, we know that $\Gamma(i) \in [-M,M]^{\binom{k}{2}}$.
Let $L$ be the least common multiple of the denominators of the (rational)
entries of $v$.  Then the entries of $L\Gamma(i)$ are integers, and
lie in the interval $[-LM, LM]$.  It follows that
$\set{ \Gamma(i) \ : \ i \geq 0}$ is a finite set of cardinality at most
$(2LM+1)^k$.  

Consider the map that sends $i$ to $\Gamma(i)$ for all $i \geq 0$.
Since this is a finite coloring of the positive integers, we know by
van der Waerden's theorem that there exist $n, d \geq 1$ such that
$\Gamma(n) = \Gamma(n+d) = \ldots = \Gamma(n + \alpha d)$.

Now $\Gamma(n+id) = \Gamma(n+ (i+1)d)$ for $0 \leq i < \alpha$, so
$$ X^{(n+id)}_l - X^{(n+id)}_m = X^{(n+(i+1)d)}_l - X^{(n+(i+1)d)}_m,$$
for $0 \leq l < m < k$
and hence
\begin{equation}
X^{(n+(i+1)d)}_l - X^{(n+id)}_l = X^{(n+(i+1)d)}_m - X^{(n+id)}_m .
\label{xn1}
\end{equation}
for $0 \leq l < m < k$.  Actually, it is easy to see that Eq.~(\ref{xn1})
holds for all $l,m $ with $0 \leq l, m < k$.

Using Eq.~(\ref{jv}), we can rewrite Eq.~(\ref{xn1}) as
$$ (\psi( {\bf w}[n+id+1..n+(i+1)d]) - dv)_l = 
 (\psi( {\bf w}[n+id+1..n+(i+1)d]) - dv)_m
$$
for $0 \leq l, m < k$.
It follows that
$$ \left| \, {\bf w}[n+id+1..n+(i+1)d] \, \right|_l - dv_l = 
 \left| \, {\bf w}[n+id+1..n+(i+1)d] \, \right|_m - dv_m $$
and hence
\begin{equation}
\left| \, {\bf w}[n+id+1..n+(i+1)d] \, \right|_l -
\left| \, {\bf w}[n+id+1..n+(i+1)d] \, \right|_m = d(v_l - v_m)
\label{wn}
\end{equation}
for $0 \leq l,  m < k$.

Now let $z = {\bf w}[n+id+1..n+(i+1)d]$.  Then Eq.~(\ref{wn}) can
be rewritten as
\begin{equation}
|z|_l - |z|_m = d(v_l - v_m)
\label{zl}
\end{equation}
for $0 \leq l, m < k$.  Note that
\begin{equation}
|z|_0 + |z|_1 + \cdots + |z|_{k-1} = |z| = d.
\label{zz1}
\end{equation}
Fixing
$l$ and summing Eq.~(\ref{zl}) over all $m \not= l$, we get 
$$ (k-1) |z|_l - \sum_{m \not= l} |z|_m = d(k-1) v_l - d \sum_{m\not= l} v_m$$
and hence by (\ref{sumv}) and (\ref{zz1}) we get
$$ (k-1) |z|_l - (d - |z|_l) = d(k-1) v_l - d(1-v_l) .$$
Simplifying, we have $k |z|_l - d = dk v_l - d$, and so $|z|_l = d v_l$.

We therefore have $ \psi({\bf w}[n+id+1..n+(i+1)d]) = d v$, for $0 \leq i < \alpha$.
Hence ${\bf w}[n+1..n+\alpha d]$ is an abelian $\alpha$-power.
\end{proof}

The following special case of Theorem~\ref{p1} is of
particular interest.

\begin{cor}
Suppose $\bf w$ is an infinite word over a finite alphabet such that in any
prefix of $\bf w$, the difference of the number of occurrences of the most
frequent letter and that of the least frequent letter is bounded by a
constant. Then $\bf w$ contains an abelian $\alpha$-power for every $\alpha
\geq 2$.
\end{cor}

\begin{proof}
Apply Theorem~\ref{p1} with $v=
({1 \over k}, {1 \over k},\ldots, {1 \over k})$.
\end{proof}

\section{Second Variation}
\label{sec3}

     Our second variation is based on the following trivial idea:  We could
avoid sum-squares if we could avoid them (mod $k$) for some integer
$k \geq 2$.
That is, instead of trying to avoid factors with blocks that sum
to the same value, we could try to avoid blocks that sum to the same value
modulo $k$.  The following result shows this is impossible, even if
we restrict our attention to blocks that sum to $0$ (mod $k$).
More general results are known 
(e.g., \cite{Justin:1972}; \cite[Chap.\ 4]{Lothaire:1983}),
but we give the proof for completeness.

\begin{theorem}
For all infinite words $\bf w$ 
over the alphabet $\Sigma_k = \lbrace 0, 1, ..., k - 1  \rbrace$
and all integers $r \geq 2$ we have that $\bf w$
contains a factor of the form $x_1 x_2 \cdots x_r$,
where $|x_1| = |x_2| = \cdots = |x_r|$ and
$\sum x_1 \equiv \sum x_2 \equiv \cdots \equiv \sum x_r \equiv
\mod{0} {k}$.
\label{just}
\end{theorem}

\begin{proof}
For $i \geq 0$ define ${\bf y}[i] = \left(
\sum_{1 \leq j \leq i} {\bf w}[i] \right) \bmod k$;
note that ${\bf y}[0] = 0$.  Then $\bf y$ is an infinite word over
the finite alphabet $\Sigma_k$, and hence by van der Waerden's theorem
there exist indices $n, n+d, \ldots, n+rd$ such that
$${\bf y}[n] = {\bf y}[n+d] = \cdots = {\bf y}[n+rd].$$
Hence ${\bf y}[n+(i+1)d] - {\bf y}[n+id] = 0$ for
$0 \leq i < r$.  But
$${\bf y}[n+(i+1)d] - {\bf y}[n+id] \equiv \mod{\sum {\bf w}[n+id+1..n+(i+1)d]}
{k},$$
so $\sum {\bf w}[n+id+1..n+(i+1)d] \equiv \mod{0} {k}$ for
$0 \leq i < r$.
\end{proof}

Theorem~\ref{just} shows that for all $k$ we cannot avoid $x x'$ 
with $|x| = |x'|$ and $\sum x \equiv \sum x' \equiv \mod{0} {k}$.
This raises the natural question, can we avoid $x x'$ 
with $|x| = |x'|$ and $\sum x \equiv \sum x' \equiv \mod{a} {k}$ for
all $a \not\equiv \mod{0} {k}$?  As phrased, the question is not so
interesting, since the word $0^\omega = 000 \cdots$ satisfies the
conditions.  If we also impose the condition
that the avoiding word be not ultimately
periodic, or even squarefree, however, then it becomes more interesting.
As we will see, we can even avoid both squares and factors
$xx'$ with $\sum x \equiv \sum x' \equiv \mod{a} {k}$ for
all $a \not\equiv \mod{0} {k}$ (with {\it no condition} on the length of
$x$ and $x'$).

\begin{theorem}
Let the morphism $\varphi$ be defined by
\begin{eqnarray*}
0 & \rightarrow & 0 \, 1 \, 0' \, {-1} \\
1 & \rightarrow & 0 \, 1 \, {-1} \, 1 \\
0' & \rightarrow & 0' \, {-1} \, 0 \, 1 \\
{-1} & \rightarrow & 0' \, {-1} \, 1 \, {-1} 
\end{eqnarray*}
and let $\tau$ be the coding defined by
\begin{eqnarray*}
0, 0' \rightarrow 0 \\
1 \rightarrow 1 \\
{-1} \rightarrow {-1} 
\end{eqnarray*}
Then the infinite word ${\bf w} = \tau(\varphi^\omega (0))$ avoids
both squares  and  factors of the form $x x'$ where
$\sum x = \sum x' \not= 0$.
\end{theorem}

\begin{proof}  The fact that $\varphi^\omega (0)$ exists follows from 
$0 \rightarrow 0 \, 1 \, 0' \, {-1}$, so that 
$\tau(\varphi^\omega (0))$ is a well-defined infinite word.

To make things a bit easier notationally,
we  may write 
$\overline{1}$ for ${-1}$.

First, let us show that $\bf w$ avoids squares.  Assume, to get a 
contradiction, that there is
such a square $xx'$ in $\bf w$, with $x = x'$,
and without loss of generality assume
$|x|$ is as small as possible.
Let $n = |x|$, and write $x = x[1..n]$, $x' = x'[1..n]$.

We call $4$ consecutive symbols of $\bf w$ that are aligned, that is, of the 
form ${\bf w}[4i+1..4i+4]$, a {\it block}.  Note that a block $B$ can be
uniquely expressed
as $\tau(\varphi(a))$ for a single symbol $a$. We call $a$ the {\it inverse image} of $B$.

Case 1:  $|xx'| \leq 25$.  It is easy to verify by exhaustive search
that all subwords of length $25$ of $\bf w$ are squarefree.
(There are only 82 such subwords.) 

\medskip

Case 2:  $|x| \geq 13$.   Then there is a block that begins at either
$x[5], x[6], x[7]$, or $x[8]$.  Such a block $y$ has at least $4$
symbols of $x$ to its left, and ends at an index at most $11$.  Thus there
are at least $2$ symbols of $x$ to the right of $y$.   We call such a
block (with at least $4$ symbols to the left, and at least $2$ to the
right) a {\it centered block}.

\medskip

Case 2a:  $|x| \equiv \mod{1,3} {4}$.   Then $x$ contains a centered block $y$.   Hence $x'$ contains an occurrence of $y$ (call it $y'$) starting at the same relative position.  Since  $|x| \equiv \mod{1,3} {4}$,
$y'$ overlaps a block $z$ starting at $1$ or $3$ positions to its left.
Since $y$ is centered, $z$ lies entirely within $x'$.  But this is impossible, since $y$ is a block, and hence starts with $0$, while the second and fourth symbol of every block $z'$ is $\pm 1$.  See Figure~\ref{fig1}.

\begin{figure}[H]
\begin{center}
\input fig1.tex
\end{center}
\caption{Case 2a}
\label{fig1}
\end{figure}

\medskip

Case 2b:  $|x| \equiv \mod{2} {4}$.  By the same reasoning, $x$ contains
a centered block $y$, so $x'$ contains an occurrence of $y$ (called $y'$)
starting at the same relative position.  Since $|x| \equiv \mod{2} {4}$,
$y'$ overlaps a block $z$ starting at $2$ positions to its left, and 
$z$ lies entirely within $x'$.  But by inspection, this can only
occur if
\begin{itemize}
\item[(i)] $y$ starts with $01$ and $z$ ends with $01$; or
\item[(ii)] $y$ starts with $0\overline{1}$ and $z$ ends with $0\overline{1}$.
\end{itemize}

In case (i), $y$ is either 
$01\overline{1}1$ or $010\overline{1}$, and $z = 0 \overline{1} 0 1$.
If $y = 0 1 \overline{1} 1$, then consider the block $z'$ that follows
$z$ in $y'$.  It must begin $\overline{1} 1$, a contradiction.  Hence
$y = 010\overline{1}$.

Now the first two symbols of $z$ precede $y'$ in $x'$ and hence must
also precede $y' = y$ in $x$.  Thus the block $y''$ that precedes $y$ in
$x$ must end in $0 \overline{1}$; it is entirely contained in $x$ because
$y$ is centered.  Hence $y'' = 0 1 0 \overline{1}$,
and $y''y$ is a shorter square in $\bf w$, a contradiction.
See Figure~\ref{fig2}.

\begin{figure}[H]
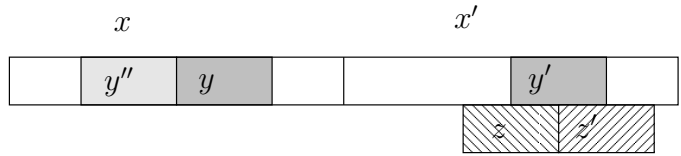

\begin{center}
\input fig2.tex
\end{center}
\caption{Case 2b(i)}
\label{fig2}
\end{figure}

In case (ii), $y$ is either
$0 \overline{1} 1 \overline{1}$ or
$0 \overline{1} 01$, and
$z = 0 1 0 \overline{1}$.  
If $y = 0 \overline{1} 1 \overline{1}$, then consider the block $z'$ that
follows $z$ in $y'$.  It must begin $1 \overline{1}$, a contradiction.
Hence $y = 0 \overline{1} 0 1$.  

Now the first two symbols of $z$ precede $y'$ in $x'$ and hence must
also precede $y'$ in $x$.  Thus the block $y''$ that precedes $y$ in
$x$ must end in $01$; it is entirely contained in $x$ because $y$ is
centered.  Hence $y'' = 0 \overline{1} 0 1$.  Hence
$y'' y$ is a shorter square in $\bf w$, a contradiction.

\medskip

Case 2c:  $|x| \equiv \mod{0} {4}$.  
Then we can write $x = r x_1 x_2 \cdots x_j l'$,
$x' = r' x'_1 x'_2 \cdots x'_j l''$, where
$lr = x_0$ (this defines $l$), $l' r' = x'_0$, $l'' r'' = x'_{j+1}$, and
$x_1, \ldots, x_j, x'_0, \ldots x'_{j+1}$ are all
blocks.  Furthermore, since $x = x'$ and
$\tau \circ \varphi$ is injective, we have
$r = r'$, $x_1 = x'_1, \ldots, x_j = x'_j$, and $l' = l''$.
See Figure~\ref{fig3}.
There are several subcases, depending
on the index $i$ in $\bf w$ in which $x$ begins.

\begin{figure}[H]
\begin{center}
\input fig3.tex
\end{center}
\caption{Case 2c}
\label{fig3}
\end{figure}

\medskip

Subcase (i):  $i \equiv \mod{1, 2} {4}$.  Then $|r| = |r'| = |r''| = 2$ or $3$.
Since any block is uniquely determined by a suffix of length $2$,
we must have $r = r'$ and so
$x_0 = x'_0$.  Hence
$x_0 \cdots x_j x_0 \cdots x_j$ corresponds to a shorter square in $\bf w$,
by taking the inverse image of each block, a contradiction.

\medskip

Subcase (ii): $i \equiv \mod{3} {4}$.  Then $|l| = |l'| = |l''| = 3$.
Again, any block is uniquely determined by a prefix of length $3$,
so $l' = l''$.  Thus $x'_0 = x'_{j+1}$ and 
$x_1 \cdots x_j x'_0 x'_1 \cdots x'_{j+1}$ is a square.  But each
of these terms is a block, so this corresponds to a shorter square 
in $\bf w$, by taking the inverse image of each block,
a contradiction.

\medskip

Subcase (iii):  $i \equiv \mod{0} {4}$.   In this case both $x$ and $x'$
can be factored into identical blocks, and hence correspond to
a shorter square in $\bf w$, by taking the inverse image
of each block, a contradiction.

\medskip
This completes the proof that $\bf w$ is squarefree.

It remains to show that if $xx'$ are consecutive factors of $\bf w$,
then $\sum x$ cannot equal $\sum x'$ unless both are $0$.

First, we prove a lemma.

\begin{lemma}
Let $\zeta$ be the morphism defined by
\begin{eqnarray*}
0, 0' & \rightarrow & 0 \, 1 \, 0' \, {-1} \\
1 & \rightarrow & 0 \, 1 \, {-1} \, 1 \, 0' \, {-1} \\
{-1} & \rightarrow & 1 \, {-1} .
\end{eqnarray*}
Then
\begin{itemize}
\item[(a)] $\varphi^n \circ \zeta = \zeta^{n+1}$ for all $n \geq 0$.

\item[(b)] $\varphi^n (0) = \zeta^n (0)$ for $n \geq 0$.
\end{itemize}
\label{lem1}
\end{lemma}

\noindent
{\it Proof.}
(a):  The claim is trivial for $n = 0$.  For $n = 1$, it becomes
$\varphi \circ \zeta = \zeta^2$, a claim that can easily be verified
by checking that $\varphi(\zeta(a)) = \zeta^2(a)$ for
all $a \in \lbrace {-1}, 0, 1, 0' \rbrace$.

Now assume the result is true for some $n \geq 1$; we prove it for $n+1$:
\begin{eqnarray*}
\varphi^{n+1} \circ \zeta &=& (\varphi \circ \varphi^n) \circ \zeta  \\
&=& \varphi \circ (\varphi^n \circ \zeta) \\
&=& \varphi \circ \zeta^{n+1}   \ \ \ \ \text{(by induction)} \\
&=& \varphi \circ (\zeta \circ \zeta^n) \\
&=& (\varphi \circ \zeta) \circ \zeta^n \\
&=& \zeta^2 \circ \zeta^n \\
&=& \zeta^{n+2}.
\end{eqnarray*}

(b):  Again, the result is trivial for $n = 0, 1$.  Assume it is true
for some $n \geq 1$; we prove it for $n+1$.    Then
\begin{eqnarray*}
\zeta^{n+1} (0) &=& \varphi^{n} (\zeta(0))  \ \ \ \ \text{(by part (a))} \\
&=& \varphi^{n} (\varphi(0)) \\
&=& \varphi^{n+1} (0).
\end{eqnarray*}
\hfill $\diamond$

\vskip 10pt
     Now let $\psi: \lbrace 0, 1, {-1} \rbrace^* \rightarrow
\lbrace 0, 1, {-1} \rbrace^*$ be defined as follows:
\begin{eqnarray*}
0 & \rightarrow & 0 \, 1 \, 0 \, {-1} \\
1 & \rightarrow & 0 \, 1 \, {-1} \, 1 \, 0 \, {-1} \\
{-1} & \rightarrow & 1 \, {-1} 
\end{eqnarray*}
Note that $\psi$ is the map obtained from $\zeta$ by equating $0$ and $0'$,
which is meaningful because $\zeta(0) = \zeta(0')$.
Then from Lemma~\ref{lem1} we get 
\begin{equation}
\tau(\varphi^n(0)) = \psi^n (0)
\label{tauv}
\end{equation}
for all $n \geq 0$.

Now form the word ${\bf v}$ from $\bf w$ by taking the
running sum.  More precisely, define ${\bf v}[i] = \sum_{0 \leq j \leq i} 
{\bf w}[j]$.  We first observe that $\bf v$ takes its values over
the alphabet $\lbrace 0, 1 \rbrace$:  From Eq.~(\ref{tauv}) we see
that ${\bf w} = \psi^\omega(0)$.  But the image of each letter under 
$\psi$ sums to $0$, and furthermore, the running sums of the image of
each letter are always either $0$ or $1$.  From this the statement about
the values of $\bf v$ follows.

Let $x x'$ be a factor of ${\bf w}$ beginning at position
$i$, with $|x| = n$, $|x'| = n'$.  
Then ${\bf w}[i..i+n-1]$ has the same sum $s$ as
${\bf w}[i+n..i+n+n'-1]$ if and only if
${\bf v}[i+n+n'] - {\bf v}[i+n] = {\bf v}[i+n] - {\bf v}[i] = s$.
In other words, 
${\bf w}[i..i+n-1]$ has the same sum $s$ as
${\bf w}[i+n..i+n+n'-1]$ if and only if
${\bf v}[i], {\bf v}[i+n],$ and ${\bf v}[i+n+n']$ form
an arithmetic progression with common difference $s$.  However,
since ${\bf v}$ takes its values in $\lbrace 0, 1 \rbrace$, this is
only possible if $s = 0$.
\end{proof}

\begin{corollary}
There exists a squarefree infinite word over $\lbrace 0,1,\ldots, k-1 \rbrace$
avoiding all factors of the form $x x'$ with
$\sum x = \sum x' = a$ for all $a \not\equiv \mod{0} {k}$.
\end{corollary}

\begin{proof}
Take the word ${\bf w} = \psi^\omega(0)$ constructed above,
and map $-1$ to $k-1$.
\end{proof}

\section{Upper and Lower Bounds}
\label{sec4}

We call a word of the form $x_1 x_2 \cdots x_r$
where  $|x_1| = |x_2| = \cdots = |x_r|$ and
$\sum x_1 \equiv \sum x_2 \equiv \cdots \equiv \mod{\sum x_r} {k}$
a {\it congruential $r$-power  (modulo $k$)}.  As we have seen,
the lengths of congruential $r$-powers, modulo $k$, are bounded.
We now consider estimating how long they can be, as a function of
$r$ and $k$.

Our first result uses some elementary number theory to get
an explicit lower bound for congruential $2$-powers.

\begin{theorem}
If $p$ is a prime, there is a word of length at least $p^2 - p - 1$ avoiding
congruential $2$-powers (modulo $p$).
\label{peng}
\end{theorem}

\begin{proof}
All arithmetic is done modulo $p$.
Let $c$ be an element of order $(p-1)/2$ in $(\Zee/(p))^*$.
If $p \equiv \mod{5,7} {8}$, let
$a$ be any quadratic residue of $p$.  If $p \equiv \mod{1,3} {8}$, let
$a$ be any quadratic non-residue of $p$.  
Let $e(k) = c^k + ak^2$ for $1 \leq k \leq p^2-p$, and define $f$ as the 
first difference of the sequence of $e$'s; that is,
$f(k) = e(k+1)-e(k)$ for $1 \leq k \leq p^2-p-1$.  Then we claim that the word
$f = f(1)f(2) \cdots f(p^2-p-1)$ avoids congruential squares (mod $p$).

To see this, assume that there is a congruential square in $f$.  Then
the sequence $e$ would have three terms where the indices and values are
both in arithmetic progression, say $k$, $k+r$, and $k+2r$.  Then
$(c^{k+r} + a(k+r)^2) - (c^k + ak^2) =
(c^{k+2r} + a(k+2r)^2) - (c^{k+r} + a(k+r)^2)$.
Simplifying, we get 
\begin{equation}
c^k (c^r - 1)^2 = -2ar^2 .
\label{eq7}
\end{equation}
If $c^r \not\equiv \mod{1} {p}$,
then 
\begin{equation}
c^k/(-2a) \equiv \mod{(r/(c^r -1))^2} {p}.
\label{eq8}
\end{equation}
Now the right-hand side of (\ref{eq8}) is a square (mod $p$), so the
left-hand side must also be a square.  But $c^k$ is a square, since
$c = g^2$ for some generator $g$.  So $-2a$ must be a square.  But 
if $p \equiv \mod{1,3} {p}$, then $-2$ is a square mod $p$, so
$-2a$ is not a square.  If $p \equiv \mod{5,7} {p}$, then $-2$ is a
nonsquare mod $p$, so $-2a$ is again not a square.  

Hence it must be that $c^r \equiv \mod{1} {p}$.  Since we chose $c = g^2$
for some generator $g$, this means that $r$ is a multiple of $(p-1)/2$,
say $r = j(p-1)/2$.  Then the left-hand side of (\ref{eq7}) is $0$ (mod $p$),
while the right hand side is $-aj^2 (p-1)^2/2$.  If this is $0$ (mod $p$),
we must have $j \equiv \mod{0} {p}$.  So $j \geq p$.  Then $2r$ is
$\geq p(p-1)$.  This gives the lower bound.
\end{proof}

\def\RA{\rightarrow}
\def\B{\hfill $\Box$}
\def\ni{\noindent}
\def\mL{\mathcal L}

We now turn to some asymptotic results.  For the remainder
of this section, as is typical in the Ramsey theory literature
\cite{Landman&Robertson:2004},
we use the language of {\it colorings}:  assigning
the $i^{\mathrm{th}}$ letter of a string $x$ to be equal to $j$ can be viewed
as coloring the integer $i$ with color $j$.

We first investigate the growth rate of the 
minimum integer $n$ such that every $k$-coloring
of $[1,n]$ admits a congruential $2$-power modulo $k$,
as $k \rightarrow \infty$.

We start with some definitions.
Let $\Omega(3,k)$ be the smallest integer $n$ such that
every set $\{x_1,x_2,\dots x_n\}$ with $x_i \in [(i-1)k+1,ik]$
contains a $3$-term arithmetic progression.
Let $\mL(k)$ be the minimum integer $n$ such that every $k$-coloring
of $[1,n]$ that uses the colors $0,1,\dots,k-1$ admits a congruential $2$-power (modulo $k$).
Finally, let
$w_1(3,k)$ be  the minimum integer $n$ such that every $2$-coloring
of $[1,n]$ admits either a $3$-term arithmetic progression
of the first color,
or $k$ consecutive  integers all with the second color.

\begin{lemma}
For any $k \in \En$, we have $\mL(k) \geq \Omega
\left(3,\big\lfloor \frac{k}{2}\big\rfloor\right)-1$.
\label{lemma1}
\end{lemma}

\begin{proof}
Consider a maximally valid set of size $n=\Omega
\left(3,\big\lfloor \frac{k}{2}\big\rfloor\right)-1$, i.e., a largest
set that avoids $3$-term arithmetic progressions.  Let
$S=\{s_1<s_2<\dots< s_n\}$ be this set.  Construct the difference set
$D=\{d_1,d_2,\dots,d_{n-1}\}=\{s_2-s_1,s_3-s_2,\dots,s_n-s_{n-1}\}$ so that $|D|=n-1$.  Note that
for any $d\in D$ we have $d \in [1,k-1]$ (so that
$0$ is not used in this construction).  We claim that
$D$ has no congruential $2$-power.  Assume, for a contradiction, that it
does.  Let $\sum_{i=x}^y d_i \equiv \sum_{y+1}^{2y-x+1} d_i \pmod k$.  Then, by
construction of $D$, we have
$$\sum_{i=x}^y d_i=s_{y+1}-s_x \quad \mbox{and}
\sum_{y+1}^{2y-x+1} d_i = s_{2y-x+2}-s_{y+1}.$$
Hence,
\begin{equation}
2s_{y+1} \equiv s_{2y-x+2}+s_x \pmod k \label{fred}. 
\end{equation}
Since $x,y+1,2y-x+2$ are in
arithmetic progression, the number of intervals
between $s_x$ and $s_{y+1}$ is the same as the
number of intervals between $s_{y+1}$ and $s_{2y-x+2}$.  Hence, 
$$
\sum_{i=x}^y d_i=s_{y+1}-s_x \in \left[(y-x)\bigg\lfloor \frac{k}{2}\bigg\rfloor+1,(y-x+2)\bigg\lfloor \frac{k}{2}\bigg\rfloor-1\right]
$$
and
$$
\sum_{y+1}^{2y-x+1} d_i = s_{2y-x+2}-s_{y+1} \in \left[(y-x)\bigg\lfloor \frac{k}{2}\bigg\rfloor+1,(y-x+2)\bigg\lfloor \frac{k}{2}\bigg\rfloor-1\right].
$$

Since the length of each of these intervals is the same and is at most
$k-1$, we see that
(\ref{fred}) is satisfied as an equality.
Hence, $s_x,s_{y+1},s_{2y-x+2}$ is a 3-term
arithmetic progression in $S$, a contradiction.
Thus, $\mL(k) > |X| = n-1 = \Omega(3,k)-2$ and we are done.
\end{proof}

Continuing, we investigate the growth rate of $\mL(k)$ through $\Omega(3,k)$.
We have the following result.

\begin{lemma}
For all $k \in \En$, $w_1(3,k) \leq k\Omega(3,k)$.
\label{lemma2}
\end{lemma}

\begin{proof}
Let $m = \Omega(3,k)$ and let $n=km$. Let $\chi$ be any
(red, blue)-coloring of $[1,n]$. Assume there are no $k$ consecutive
blue integers. So, for each $i$,  $1 \leq i \leq m$, the interval
$[(i-1)k+1,ik]$ contains a red element, say $a_i$. Then, by the
definition of $\Omega(3,k)$, there is a $3$-term arithmetic progression among
the
$a_i$'s.
\end{proof}

Recently, Ron Graham \cite{Graham:2006} has shown the following.

\begin{theorem}
\rm (Graham) \it
There exists a constant $c>0$
such that, for $k$ sufficiently large,
$w_1(3,k) >k^{c\log k}$.
\label{thm3}
\end{theorem}

As a corollary, using Lemma~\ref{lemma2}, we have

\begin{corollary}
There exists a constant $c>0$
such that, for $k$ sufficiently large,
$\Omega(3,k) >k^{c\log k}$.
\label{corollary4}
\end{corollary}

\begin{proof}
From Theorem~\ref{thm3} and Lemma~\ref{lemma2} we have, for some $d>0$,
$$
\Omega(3,k)\geq \frac{w_1(3,k)}{k} > k^{d\log k - 1} > k^{\frac{d}{2}\log k}.
$$
Taking $c=\frac{d}{2}$ gives the result.
\end{proof}

We now apply Corollary~\ref{corollary4} to Lemma~\ref{lemma1}
to yield the following theorem, which
states that $\mL(k)$ grows faster than any polynomial in $k$.

\begin{theorem}
There exists a constant $c>0$
such that, for $k$ sufficiently large,
$\mL(k) >k^{c\log k}$.
\label{theorem5}
\end{theorem}

\begin{proof}
We have (suppressing constant terms)
$$
\mL(k) \geq \Omega\left(3,\bigg\lfloor \frac{k}{2}\bigg\rfloor\right)
> \left(\frac{k}{2}\right)^{d \log \frac{k}{2}}
$$
for some $d>0$, provided $k$ is sufficiently large.  Since
$\frac{k}{2} > \sqrt{k}$ for $k>4$ this gives, for sufficiently large $k$,
$$
\mL(k) > k^{\frac{d}{2} \log \frac{k}{2}} > k^{\frac{d}{4} \log k}.
$$
Taking $c=\frac{d}{4}$ yields the result.
\end{proof}

We now turn from congruential $2$-powers to the more general case
of congruential $t$-powers.
To this end, define
$\mL(k,t)$ to be the minimum integer $n$ such that every $k$-coloring
of $[1,n]$ using the colors $0,1,\dots,k-1$
admits a congruential $t$-power modulo $k$.

Adapting the proof of Lemma~\ref{lemma1}, to this we immediately get

\begin{lemma}
For any $k,t \in \En$, we have $\mL(k,t) \geq \Omega
\left(t + 1,\big\lfloor \frac{k}{2}\big\rfloor\right)-1$.
\label{lemma6}
\end{lemma}

Now, a result due to Nathanson \cite{Nathanson:1980} gives us the following result.

\begin{theorem}
For any $k,t \in \mathbb{Z}^+$, we have
$$
\Omega
\left(t + 1,\Big\lfloor \frac{k}{2}\Big\rfloor\right) \geq 
w \left( \Big\lceil \frac{2t}{k}\Big\rceil+1; \Big\lfloor \frac{k}{2}\Big\rfloor\right).
$$
\label{theorem7}
\end{theorem}

When $k = 4$, this gives us the following.

\begin{corollary}
For any $t \in \mathbb{Z}^+$ we have $\mL(4,t) \geq 
w\left( \Big\lceil \frac{t}{2}\Big\rceil+1; 2\right) - 1$.
\label{corollary8}
\end{corollary}

Hence, this says, roughly, that $\mL(4,2\ell)$ serves as an upper bound
for the classical van der Waerden number $w(\ell,\ell)$.

A recent result of Bourgain \cite{Bourgain:2008} implies the bound
$w(3;k) = o(k^{ck^{3/2}})$ for some constant $c>0$.

Hence, for sufficiently large $k$, there exist constants $c,d>0$
such that
$$k^{c \log k} < \mL(k) < k^{d k^{3/2}}$$
so that we have a very rough idea of the growth rate. 

\section{Computational Results}
\label{sec5}

As we have seen,
the known upper bounds on van der Waerden numbers provide upper
bounds for the length of the longest word avoiding congruential powers.
We also did some explicit computations.
We computed the length $l(r,k)$ of the longest word over $\Sigma_k$ avoiding
congruential $r$-powers (modulo $k$),
for some small values of $k$ and $r$, and the
lexicographically least such longest word $x_{r,k}$. The data are summarized below.

\begin{center}
\begin{tabular}{|l|l|r|l|}
\hline
$r$ & $k$ & $l(r,k)$ & $x_{r,k}$  \\
\hline
2 & 2 &  3 & 010 \\
2 & 3 &  7 & 0102010 \\
2 & 4 & 16 & 0130102013101201 \\
2 & 5 & 33 & 010214243213143040102142432131430 \\
2 & 6 & 35 & 01024021240241402401024021240241402 \\
2 & 7 & 47 & 01021614636032312426404301021614636032312426404 \\
\hline
3 & 2 &  9 & 001101100 \\
3 & 3 & 67 & {\tiny 0010210112021200102022121011202120010201012101120212001021002210112} \\
\hline
4 & 2 & 88 & {\tiny 0011000110001001110010001100011000100111001000110001100010011100100011000110001001110011}\\
\hline
\end{tabular}
\end{center}

It remains an interesting open problem to find better upper and lower
bounds on the length of the longest word avoiding congruential powers.

\section{Acknowledgments}
\label{sec6}
We thank Thomas Stoll for a helpful observation.  Theorem~\ref{peng} was inspired by some related empirical calculations
by Richard (Yang) Peng.

\end{document}